\documentclass[11pt,a4paper,leqno]{article}
\usepackage[latin1]{inputenc} 
\usepackage{amsmath}
\usepackage{amsfonts}
\usepackage{amssymb}
\usepackage{MnSymbol}
\usepackage{stmaryrd}

\newtheorem{satz}{Satz}[section]
\newtheorem{theorem}[satz]{Theorem}
\newtheorem{lemma}[satz]{Lemma}
\newtheorem{prop}[satz]{Proposition}
\newtheorem{cor}[satz]{Corollary}

\newcommand{\rund}[1]{\left(#1\right)}
\newcommand{\spitz}[1]{\left\langle{#1}\right\rangle}

\newcommand{\schweif}[1]{\left\{#1\right\}}

\def\cz{\mathbb{C}}
\def\nz{{\rm I\kern-.20em N}}
\def\rz{{\rm I\kern-.20em R}}

\def\M{\mathcal{M}}

\def\O{\mathcal{O}}
\def\T{\mathcal{T}}
\def\S{\mathcal{S}}

\def\m{\mathfrak m}

\def\Diff{{\rm Diff}}
\def\Mod{{\rm Mod}}

\def\Aut{{\rm Aut}}

\def\id{{\rm id}}

\def\rel{{\rm rel}}

\def\strong{{\rm strong}}

\def\eps{\varepsilon}

\begin{document}

\begin{center}

{\Large \bf Preliminaries to}

\vskip 0.5 true cm

{\Huge \bf Versal Families of Compact Super

\vskip 0.4 true cm

Riemann Surfaces}

\vskip 1.0 true cm

Roland Knevel \\

Bar-Ilan University RAMAT GAN

Department of Mathematics 

ISRAEL
\end{center}

\vskip 0.8 true cm

We collect some classical results on versal families of ordinary compact Riemann surfaces needed in \cite{Knevel}.

\vskip 0.8 true cm

{\bf Grauert's theorem:} $W$ Stein manifold, $E, F \twoheadrightarrow W$ holomorphic vector bundles. Then $E \simeq F$ as holomorphic vector bundles iff $E \simeq F$ as topological vector bundles. In particular all holomorphic vector bundles on a contractible Stein manifold are trivial. \\

see \cite{Ballico} introduction. \\

%{\bf Cartan's theorem A:} $\S$ coherent sheaf on a Stein manifold $W$~. Then for all $w \in W$~, $\S_w$ as $\rund{\O_W}_w$-module is generated by $\schweif{\eckig{f}_w \, | \, f \in \S(W)}$~. \\

{\bf Cartan's theorem B:} $\S$ coherent sheaf on a Stein manifold $W$~. Then $H^k\rund{W, \S} = 0$ for all $k \geq 1$~. \\

see \cite{GrauRemStein} introduction. \\

{\bf \cite{GrauRem} theorem 10.5.5:} If $\dim_\cz H^i\rund{X_y, \underline V_y}$ is independent of $y \in Y$ then all sheaves $f_{(i)}\rund{\underline V}$ are locally free and all maps 

\[
f_{y, i}: f_{(i)}\rund{\underline V} \left/ \m_y f_{(i)}\rund{\underline V}\right. \rightarrow H^i\rund{X_y, \underline V_y}
\]

are isomorphisms. \\

Hereby $f: X \twoheadrightarrow Y$ is a holomorphic family of compact complex manifolds $X_y := f^{- 1}(y)$~, $y \in Y$~, and $\underline V \twoheadrightarrow X$ a holomorphic vector bundle. $\underline V_y := \left. \underline V \right|_{X_y}$~, and $\m_y \lhd \O_Y$ denotes the maximal ideal of all holomorphic functions vanishing at the point $y \in Y$~. Finally $f_{y, i}: f_{(i)}\rund{\underline V} \left/ \m_y f_{(i)}\rund{\underline V}\right. \rightarrow H^i\rund{X_y, \underline V_y}$ denotes the canonical homomorphism.

\begin{lemma} \label{vanishH1} $W$ simply connected manifold, $G$ group. Then $H^1(W, G)$ is trivial. \end{lemma}

{\it Proof:} Let $W = \bigcup_i U_i$ be an open cover and $g_{i j}: U_i \cap U_j \rightarrow G$ locally constant, $U_i \cap U_j \not= \emptyset$~, forming a $1$-cocycle, so $g_{k i} g_{j k} g_{i j} = 1$ if $U_i \cap U_j \cap U_k \not= \emptyset$~, and $g_{i i} = 1$~. We have to show that after maybe refining the open cover there exist $h_i: U_i \rightarrow G$ locally constant such that $g_{i j} = h_j^{- 1} h_i$~.

After maybe refining the open cover we may assume by \cite{BottTu} corollary 5.2 that it is good, which means that all finite non-empty intersections of the $U_i$ are diffeomorphic to $\rz^n$~. Then of course all $g_{i j}$ are constant.

Choose $i_0$~. Take an arbitrary $i$ and a path $\gamma: [0, 1] \rightarrow W$ from $U_{i_0}$ to $U_i$~. Let $I_{k \sigma}$~, $\sigma = 1, \dots, s_k$~, numbered increasingly, be the connected components of $\gamma^{- 1}\rund{U_k}$~. Then $[0, 1] = \bigcup_{k, \sigma} I_{k \sigma}$ is a cover with open intervals, and there exist unique $h_{k \sigma} \in G$~, such that $h_{i_0, 1} = 1$ and $g_{k l} = h_{l \tau}^{- 1} h_{k \sigma}$~, $I_{k \sigma} \cap I_{l \tau} \not= \emptyset$~. Since two such paths are homotopic and $h_{i, s_i}$ is constant under small variation of $\gamma$~, we see that infact $h_i := h_{i, s_i}$ is independent $\gamma$~.

$\leadsto$ $\rund{h_i}$ 0-chain in $G$~, and taking a path from $U_{i_0}$ to $U_j$ over $U_i \cap U_j$ one sees that indeed $g_{i j} = h_j^{- 1} h_i$~. $\Box$

\begin{theorem} \label{smoothfamily} $\pi: M \twoheadrightarrow \rz^n$ smooth family of compact real manifolds. Then $M$ is trivial. \end{theorem}

{\it Proof:} by induction on $n$, inspired by \cite{Kodaira} proof of theorem 2.4.

$n = 0$ trivial.

$n \rightarrow n + 1$~: We construct a smooth vectorfield $\chi \in H^0\rund{T M}$ such that $(d \pi) \chi = \partial_1 \circ \pi$~: indeed locally in $M$ such $\chi$ exist, so we obtain an open cover $M = \bigcup_i U_i$~, without restriction locally finite, and $\chi_i \in H^0\rund{T U_i}$ with this property. Now take a smooth partition $\rund{\eps_i}$ of unity subordinated to $\rund{U_i}$ and $\chi := \sum_i \eps_i \chi_i$~.

Let $\Phi$ be the integral flow to $\chi$~, so the maximal solution of the initial value problem $\Phi|_{\{0\} \times M} = \id_M$ and $(d \Phi) \partial_u = \chi \circ \Phi$~. Since $(d \pi) \chi = \partial_1 \circ \pi$ we obtain

\[
\begin{array}{ccc}
\phantom{12345678} \rz \times M & \mathop{\longrightarrow}\limits^\Phi & M \phantom{12345678901} \\
\rund{\id_\rz, \pi} \twoheaddownarrow & \circlearrowleft & \twoheaddownarrow \pi \phantom{123456789,} \\
\phantom{12345678} \rz \times \rz^{n + 1} & \longrightarrow & \rz^{n + 1} \phantom{1234567890,} \\
\phantom{123456789} (u, t) & \mapsto & t + (u, 0, \dots, 0)
\end{array} \,~.
\]

So since all fibers of $M$ are compact we see that $\Phi$ is defined on all $\rz \times M$~. Let $\pi_1: M \twoheadrightarrow \rz$ denote the first component of $\pi$~. Then

\[
\Phi|_{\rz \times M|_{\{0\} \times \rz^n}}: \rz \times M|_{\{0\} \times \rz^n} \rightarrow M
\]

is a strong family morphism with inverse $x \mapsto \rund{\pi_1(x), \Phi\rund{- \pi_1(x), x}}$~, so a strong family isomorphism. Finally by induction hypothesis $M|_{\{0\} \times \rz^n}$ is trivial. $\Box$

\begin{theorem} \label{autHomotopic} $X$ compact Riemann surface of genus $g \geq 2$~, $\Phi \in \Aut \ X$ homotopic to $\id_X$~. Then $\Phi = \id_X$~. \end{theorem}

{\it Proof:} see \cite{Hurwitz}. $\Box$

% Without restriction $X = X_{u_0}$~, $u_0 \in \T_g$~. $\Phi$ can be extended to $(\Gamma, \gamma) \in \Aut \ M_g$ by \cite{Grothendieck} proposition 3.3. Since $\Phi$ is homotopic to $\id_{X_{u_0}}$ also $\Gamma|_{X_u}: X_u \mathop{\rightarrow}\limits^\sim X_{\gamma(u)}$ is homotopic to $\id_{X_{u_0}}$ and therefore $\gamma(u) = u$ for all $u \in \T_g$ by construction of $\T_g$~, which means that $\Gamma$ is strong. By theorem \ref{autStrong} if $g \geq 3$ then $\Gamma = \id_{M_g}$~, if $g = 2$ then $\Gamma = \id_{M_g}$ or $\Gamma = J$~, but $J$ is not homotopic to $\id_{M_2}$~, see \cite{Ares} section 2.3. $\Box$

\begin{cor} \label{ideverywhere} $N \twoheadrightarrow W$ holomorphic family of compact Riemann surfaces $Y_w$ of genus $g \geq 2$~, $W$ connected, $\Phi \in \Aut_\strong \ N$~, $w_0 \in W$ such that $\Phi|_{Y_{w_0}} = \id_{Y_{w_0}}$~. Then $\Phi = \id_N$~. \end{cor}

{\it Proof:} Since $W$ is connected all $\Phi|_{Y_w}$~, $w \in W$~, are homotopic to $\id_{Y_w}$~, so equal to $\id_{Y_w}$ by theorem \ref{autHomotopic}. $\Box$ \\

% Take a local chart $\pr_{\Omega}: \Omega \times U \twoheadrightarrow \Omega$~, $\Omega \subset \cz^n$~, $U \subset \cz$ open, $w_0 \in \Omega$~, of $N \twoheadrightarrow W$~. In this chart

%\[
%\Phi = \rund{w, z + \sum_{n = 1}^\infty w^n \varphi_n(z)} \,~,
%\]

%locally in $W$ around $w_0$~, $z$ the coordinate on $\cz$~, $w$ the coordinates on $\cz^n$~, all $\varphi_n$ holomorphic. Infact all $\varphi_n = 0$ by induction on $n$~:

%\begin{quote}
%Assume $\varphi_\nu = 0$~, $\nu < n$~, so $\Phi = \rund{w, z + \sum_{\nu = n}^\infty w^\nu \varphi_\nu(z)}$~. By induction on $k$ :

%$\Phi^k = \rund{w, z + k w^n \varphi_n(z) + \text{ higher terms in } w}$~. On the other hand $\abs{\Aut \ X} \leq~84~(g - 1)$ for all compact Riemann surfaces $X$ of genus $g \geq 2$ by Hurwitz' estimate, see \cite{FarkKra} theorem V.1.3. We see that $\Phi^{(84 (g - 1))!} = \id_N$~, and so $(84 (g - 1))! \ \varphi_n = 0$~.
%\end{quote}

%By the identity theorem $\Phi = \id_N$~. $\Box$ \\

For every $g \in \nz$ let $\T_g$ denote the Teichmüller space for genus $g$~, constructed as set of equivalence classes of marked compact Riemann surfaces of genus $g$~, see for example \cite{Ares} section 2.3.

\begin{prop}
\item[(i)] For every genus $g$ there exists a holomorphic family $M_g \twoheadrightarrow \T_g$ of compact Riemann surfaces $X_u$ such that the class of $X_u$ equals the point $\overline u$ in $\Mod_g$~.
\item[(ii)] $\T_g$ is a bounded contractible domain of holomorphy, so in particular Stein.
\item[(iii)] $M_g$ is trivial as a smooth family of compact smooth families.
\end{prop}

{\it Proof:} (i) trivial for $g = 0$~: $\T_0$ single point and $M_0 = \widehat \cz$~. Explicit construction for $g = 1$~: $\T_1$ upper half plane, $M_1 = \left.\T_1 \times \cz \right/ \spitz{S, T}$~, $S: (u, z) \mapsto (u, z + 1)$~, $T: (u, z) \mapsto (u, z + u)$~. See \cite{Duma} section 1 for $g \geq 2$~.

(ii) trivial for $g = 0$ and $g = 1$~. For $g \geq 2$ by \cite{Earle} section 2.

(iii) trivial for $g = 0$ and $g = 1$~. For $g \geq 2$ by theorem \ref{smoothfamily} since $\T_g$ diffeomorphic to $\rz^{6 (g - 1)}$~, see \cite{Natanzon} theorem 4.1. $\Box$

\begin{theorem}[Anchoring property of $M_g$ ] \label{anchoring} Let $\pi_N: N \twoheadrightarrow W$ be a holomorphic family of compact Riemann surfaces $Y_w$ of genus $g$~, $W$ simply connected if $g \geq 2$~, $W$ contractible and Stein if $g \leq 1$~, $\sigma: Y_{w_0} \mathop{\rightarrow}\limits^\sim X_{u_0}$~, $w_0 \in W$ and $u_0 \in \T_g$~. Then there exists a unique $\varphi: W \rightarrow \T_g$ such that $\varphi\rund{w_0} = u_0$ and $\varphi$ can be extended to a fiberwise biholomorphic family morphism $(\Phi, \varphi): N \rightarrow M_g$ with $\Phi|_{Y_{w_0}} = \sigma$~. If $g \geq 2$ then also $\Phi$ is uniquely determined by~$\sigma$~. \end{theorem}

{\it Proof:} By theorem \ref{smoothfamily} there exists an open cover $W = \bigcup_i \Omega_i$ and trivializations \\
$\varphi_i: N|_{\Omega_i} \mathop{\rightarrow}\limits^\sim \Omega_i \times Y_{w_0}$ as smooth families of compact smooth manifolds. On $\Omega_i \cap \Omega_j$~, $\varphi_i$ and $\varphi_j$ differ by the strong smooth family automorphism $\varphi_j \circ \varphi_i^{- 1}$ of $\rund{\Omega_i \cap \Omega_j} \times Y_{w_0}$~, and the homotopy class of $\rund{\varphi_j \circ \varphi_i^{- 1}}(w, \diamondsuit) \in \Diff \ Y_{w_0}$ is locally independent of $w$~. $W$ is simply connected, so $H^1\rund{W, \rund{\Diff \ Y_{w_0}} \left/ \rund{\Diff_0 \ Y_{w_0}}\right.} = 0$ by lemma \ref{vanishH1}, where $\Diff_0 \ Y_{w_0}$ denotes the normal subgroup consisting of diffeomorphisms homotopic to $\id_{Y_{w_0}}$~. So after refining the open cover $\rund{\Omega_i}$ we may assume without restriction that all $\rund{\varphi_j \circ \varphi_i^{- 1}}(w, \diamondsuit) \in \Diff_0 \ Y_{w_0}$~. Choose $i_0$ with $w_0 \in \Omega_{i_0}$~. Then by taking $\rund{\id_{\Omega_{i_0}}, \rund{\left.\varphi_{i_0}\right|_{Y_{w_0}} }^{- 1}} \circ \varphi_i$ instead of $\varphi_i$ we ensure that $\left.\varphi_i\right|_{Y_{w_0}} \in \Diff_0 \ Y_{w_0}$ whenever $w_0 \in \Omega_i$~. From now on we identify $N|_{\Omega_i}$ with $\Omega_i \times Y_{w_0}$ via $\varphi_i$ as smooth families, and for all $w \in W$ the homotopy class of $\left.\varphi_i\right|_{Y_w}: Y_w \rightarrow Y_{w_0}$ is independent of the choice of $i$~.

By construction of $\T_g$~, for every $w \in W$ there exists a unique $u_w \in \T_g$ such that there exists an isomorphism $Y_w \mathop{\rightarrow}\limits^\sim X_{u_w}$ homotopic to~$\sigma$~. Of course $u_{w_0} = u_0$~. \\

{\it Uniqueness:} Take $(\Phi, \varphi)$ as desired. Then $\Phi|_{Y_w}: Y_w \mathop{\rightarrow}\limits^\sim X_{\varphi(w)}$ is homotopic to $\sigma$~, so $\varphi(w) = u_w$ for all $w \in W$~. If $g \geq 2$ then also $\Phi$ is unique by corollary \ref{ideverywhere}. \\

{\it Existence:} $g \geq 2$ : By \cite{Duma} section 1 there exist, locally in $W$~, fiberwise biholomorphic family morphisms $(\Phi, \varphi): N \rightarrow M_g$~. Also by \cite{Duma} section 1 every orientation preserving homeomorphism of $Y_{w_0}$ induces a family automorphism of $M_g$~, so without restriction $\Phi|_{Y_w}: Y_w \mathop{\rightarrow}\limits^\sim X_{\varphi(w)}$ is homotopic to $\sigma$~, which implies $\varphi(w) = u_w$ by construction of $\T_g$~. Therefore $\varphi$ is infact globally defined. Now on their overlaps two local choices of $\Phi$ differ by a strong automorphism of $\varphi^* M_g$ homotopic to $\id_{\varphi^* M_g}$~, which so equals $\id_{\varphi^* M_g}$ by theorem \ref{autHomotopic}. We see that also $\Phi$ is globally defined. Finally $\Phi|_{Y_{w_0}} \circ \sigma^{- 1} \in \Aut \ X_{u_0}$ homotopic to $\id_{X_u}$~, so $\Phi|_{Y_{w_0}} = \sigma$ by theorem \ref{autHomotopic}. \\

$g = 1$ : $\rund{\pi_N}_* \rund{T^\rel N}^*$ is globally free of rank $1$ by \cite{GrauRem} theorem 10.5.5 and Grauert's theorem since $W$ is contractible and Stein, so let $\omega$ be a global frame such that $\omega\rund{w_0, \diamondsuit} = \sigma^* d z$~. Let $A$ and $B$ be the pullbacks of the straightlines $[0, 1]$ resp. $[0, u_0]$ under $\sigma \circ\left.\varphi_i\right|_{Y_w}$~. They are closed curves in $Y_w$~, whose homotopy classes are independent of $i$ and generate $\pi_1\rund{Y_w}$~. Therefore

\[
\psi, \varphi: W \rightarrow \cz \,~, \, w \mapsto \int_A \omega(w, \diamondsuit) \text{ resp. } \int_B \omega(w, \diamondsuit)
\]

are well-defined and holomorphic on $W$~, $\psi(w), \varphi(w)$ linearly independent over $\rz$ for all $w \in W$~. $\psi\rund{w_0} = 1$ and $\varphi\rund{w_0} = u_0$~. After multiplying $\omega$ with $\frac{1}{\psi}$ we may assume that $\psi = 1$~. Then $\varphi(W) \subset \T_1$~. Locally in $W$ take a holomorphic section $q$ of $N$~. So locally in $W$~, $(\Phi, \varphi)$ given by

\[
z \mapsto \rund{\varphi\rund{\pi_N(z)}, \int_q^z \omega\rund{\pi_N(z), \diamondsuit}}
\]

is a fiberwise biholomorphic family morphism $N \rightarrow M_1$~. On their overlaps two such locally constructed $\Phi$ differ by a translation with a holomorphic section of $\varphi^* M_1$~. Since $W$ contractible and Stein, $H^1\rund{W, \varphi^* M_1} = 0$ by Cartan's theorem B, and so $\Phi$ can be defined globally. Finally $\Phi|_{Y_{w_0}} \circ \sigma^{- 1} \in \Aut \ X_{u_0}$ is a translation $\tau_a$ with some $a \in \cz$~, and so taking $\tau_{- a} \circ \Phi$ instead of $\Phi$ $\, \leadsto$ $\Phi|_{Y_{w_0}} = \sigma$~. \\

$g = 0$ : Locally in $W$ take holomorphic line bundles $L \twoheadrightarrow N$ such that $L^{\otimes 2} = T^\rel N$~. $\leadsto$ open cover $W = \bigcup_i U_i$ and holomorphic line bundles $L_i \twoheadrightarrow N|_{U_i}$ such that $L_i^{\otimes 2} = T^\rel N$~. After refining the open cover we may assume that it is good, see \cite{BottTu} corollary 5.2.

$\left.\rund{L_i \otimes L_j^*}\right|_{Y_w}$ is trivial for all $w \in U_i \cap U_j$~. Therefore $\rund{\pi_N}_* \rund{L_i \otimes L_j^*}$ is a locally free sheaf on $U_i \cap U_j$ of rank $1$ by \cite{GrauRem} theorem 10.5.5, so every local frame yields an isomorphism \\
$\psi_{i j}: L_i \mathop{\rightarrow}\limits^\sim L_j$ locally in $U_i \cap U_j$~, and we may assume that $\psi_{i j}^{\otimes 2} = \id_{T^\rel N}$~. Infact we can construct a global $\psi_{i j}$ on $U_i \cap U_j$~: Two local choices of $\psi_{i j}$ differ by multiplication with $\pm 1$~, and so we see that the obstructions to define $\psi_{i j}$ globally on $U_i \cap U_j$ lie in $H^1\rund{U_i \cap U_j, \schweif{\pm 1}}$~, which vanishes since $U_i \cap U_j$ is contractible.

Without restriction $\psi_{j i} = \psi_{i j}^{- 1}$~. Then

\[
\varphi_{i j k} := \psi_{k i} \circ \psi_{j k} \circ \psi_{i j} \in \Aut \ L_i
\]

form a $2$-cocycle in $\schweif{\pm 1}$~. $H^2\rund{W, \schweif{\pm 1}} = 0$ since $W$ is contractible, and so without restriction all $\varphi_{i j k} = 1$~. Therefore all $L_i$ glue together to a holomorphic line bundle $L \twoheadrightarrow N$ with $L^{\otimes 2} = T^\rel N$~, so in particular $L$ is of rank $1$~.

By \cite{GrauRem} theorem 10.5.5 and Grauert's theorem since $W$ is contractible and Stein $\rund{\pi_N}_* L \simeq~\O_W^{\oplus 2}$~, so let $(f, h)$ be a global frame. $\Phi := \frac{f}{h} \in \M(N)$~, $\Phi|_{Y_w}$ is holomophic apart from one single pole for every $w$~. So $\Phi: N \rightarrow \widehat\cz$ is fiberwise biholomorphic. Com\-po\-sing $\Phi$ with a suitable element of $\Aut \ \widehat\cz = PSL(2, \cz)$ $\, \leadsto$ $\Phi|_{Y_{w_0}} = \sigma$~. $\Box$

\begin{theorem}[Infinitesimal universality of $M_g$ ] $\rund{d M_g}_u: T_u \T_g \rightarrow H^1\rund{T X_u}$ is an isomorphism for all $u \in \T_g$~. \end{theorem}

{\it Proof:} {\it Surjectivity:} Let $\beta \in H^1\rund{T X_u}$~. Since $H^2\rund{T X_u} = 0$~, by \cite{Kodaira} theorem 5.6 there exists a holomorphic family $N \twoheadrightarrow \Delta$ of compact Riemann surfaces $Y_w$~, $\Delta \subset \cz$ a disc around $0$~, such that $Y_0 = X_u$ and $\beta = (d N)_0 \partial_w$~. By the anchoring property of $M_g$ there exists a fiberwise biholomorphic family morphism $(\varphi, \Phi): N \rightarrow M_g$ with $\varphi(0) = u$~, $\Phi|_{Y_0} = \id_{Y_0}$~, and so by the chain rule $\beta = (d N)_0 \partial_w = \rund{d M_g}_u (d \varphi)_0 \partial_w$~.

{\it Injectivity:} by equality of dimension. $\Box$

\begin{theorem} The modular group $\Gamma_g$ acts properly discontinuously on $\T_g$~. \end{theorem}

{\it Proof:} obvious for $g \leq 1$~, see \cite{Duma} section 2 for $g \geq 2$~. $\Box$
\begin{theorem} \label{autStrong} $\schweif{\id_{M_1}, J}$ is a set of representatives for $\left.\rund{\Aut_\strong \ M_1} \right/ \schweif{\text{translations}}$~, \\
$J: (w, z) \mapsto (w, - z)$~,

$\Aut_\strong \ M_2 = \schweif{\id_{M_2}, J}$~, $J$ the hyperelliptic reflection, and

$\Aut_\strong \ M_g = \schweif{\id_{M_g}}$ for $g \geq 3$~.
\end{theorem}

{\it Proof:} simple calculation for $M_1$~.

By \cite{Ares} section 2.3 the hyperelliptic reflection $J$ is in $\Aut_\strong \ M_2 \setminus \schweif{\id_{M_2}}$~. Let \\
$\Phi \in \Aut_\strong \ M_2$~. By \cite{Poonen} theorem 1 there exists $u \in \T_2$ such that $\Aut \ X_u = \schweif{\id_{X_u}, J|_{X_u}}$~, and so $\Phi = \id_{M_2}$ or $\Phi = J$ by corollary \ref{ideverywhere}.

Finally let $\Phi \in \Aut_\strong \ M_g$~, $g \geq 3$~. By \cite{Baily} there exists $u \in \T_g$ such that $\Aut \ X_u = \schweif{\id_{X_u}}$~, and so $\Phi = \id_{M_g}$ by corollary \ref{ideverywhere}. $\Box$


\begin{thebibliography}{99}

\bibitem{Ares} {\sc Arés}, P.: Notes on Teichmüller spaces

http://www.math.tifr.res.in/\texttildelow pablo/download.php~.

\bibitem{Baily} {\sc Baily}, W. L. Jr.: On the automorphism group of a generic curve of genus $> 2$~. J. Math. Kyoto Univ. {\bf 1}-1 (1961), 101 - 108.

\bibitem{Ballico} {\sc Ballico}, E.: Holomorphic vector bundles on holomorphically convex complex manifolds. Georgian Math. J. {\bf 13}-1 (2006), 7 - 10.

\bibitem{BottTu} {\sc Bott}, R. and {\sc Tu} L. W.: Differential Forms in Algebraic Topology, Springer, New York, 1982.

\bibitem{Duma} {\sc Duma}, A.: Die Automorphismengruppe der universellen Familie kompakter Riemannscher Flächen vom Geschlecht $g \geq 3$~. manuscripta math. {\bf 17} (1975), 309 - 315.

\bibitem{Earle} {\sc Earle}, C. J.: Families of Riemann surfaces and Jacobi varieties. Ann. of Math., {\bf 107} (1978), 255 - 286.

\bibitem{FarkKra} {\sc Farkas}, H. M. and {\sc Kra}, I.: Riemann Surfaces. Springer New York 1991.

\bibitem{GrauRem} {\sc Grauert}, H. and {\sc Remmert}, R.: Coherent Analytic Sheaves. A Series of Comprehensive Studies in Mathematics. Springer Berlin Heidelberg 1984.

\bibitem{GrauRemStein} {\sc Grauert}, H. and {\sc Remmert}, R.: Theory of Stein spaces. A Series of Comprehensive Studies in Mathematics, Springer, Berlin Heidelberg, 1977.

%\bibitem{Grothendieck} {\sc Grothendieck}, A.: Techniques de construction en géométrie analytique. I. Description axiomatique de l'espace de Teichmüller et de ses variantes. Séminaire Henri Cartan, tome 13, $\text{n}^\text{o}$ 1 (1960 - 1961), exp $\text{n}^\text{o}$ 7 et 8, 1 - 33.

\bibitem{Hurwitz} {\sc Hurwitz}, A.: Über algebraische Gebilde mit eindeutigen Transformationen in sich. Math. Ann. {\bf 41}(1893), 403-442.

\bibitem{Knevel} {\sc Knevel}, R.: Versal Families of Compact Super Riemann Surfaces. Preprint 2012.

\bibitem{Kodaira} {\sc Kodaira}, K.: Complex manifolds and deformation of complex structures. Springer Heidelberg 2005.

\bibitem{Natanzon} {\sc Natanzon}, S.: Moduli of Riemann surfaces, Hurwitz-type spaces, and their superanalogues. Russian Math. Surveys {\bf 54}-1 (1999), 61 - 117.

\bibitem{Poonen} {\sc Poonen}, B.: Varieties without extra isomorphisms II: Hyperelliptic curves. Math. Res. Lett. {\bf 7}-1 (2000), 77-82.

\end{thebibliography}
\end{document}